\documentclass{amsart}
\usepackage{amssymb}
\usepackage{amsmath}

\newtheorem{theorem}{Theorem}
\newtheorem{Main theorem}{Main theorem}

\newtheorem{lemma}{Lemma}
\newtheorem{corollary}{Corollary}

\newcommand{\LL}{\ensuremath{\mathbb{L}}}

\newcommand{\N}{\ensuremath{\mathbb{N}}}
\newcommand{\Z}{\ensuremath{\mathbb{Z}}}
\newcommand{\Q}{\ensuremath{\mathbb{Q}}}

\newcommand{\Pro}{\ensuremath{\mathbb{P}}}
\newcommand{\A}{\ensuremath{\mathbb{A}}}

\begin{document}
\title{Arcs and resolution of singularities}
\author{Johannes Nicaise$^\dag$}

\email{johannes.nicaise@wis.kuleuven.ac.be}
\address{Department of Mathematics\\
Katholieke Universiteit Leuven\\ Celestijnenlaan 200B\\ B-3001
Leuven\\ Belgium}
\thanks{$\dag$Research Assistant of the Fund for Scientific Research --
 Flanders (Belgium)(F.W.O.)
}
\begin{abstract}
For a certain class of varieties $X$, we derive a formula for the
valuation $d_{X}$ on the arc space $\mathcal{L}(Y)$ of a smooth
ambient space $Y$, in terms of an embedded resolution of
singularities. A simple transformation rule yields a formula for
the geometric Poincar\'e series. Furthermore, we prove that for
this class of varieties, the arithmetic and the geometric
Poincar\'e series coincide. We also study the geometric valuation
for plane curves.

\end{abstract}
\maketitle
\section{Introduction}
Let $k$ be a field of characteristic zero, and let $k^{alg}$ be an
algebraic closure. Let $X$ be a subvariety of affine space
$\A^{d}_{k}$, defined by polynomial equations
$f_j(x)=0,\,j=1,\ldots ,r$.

We can describe jets on $X$ in terms of the coordinate system
$(x_1,\ldots,x_d)$ on $\A^{d}_k$. An $n$-jet on $X$ is a tuple of
truncated power series $$a=(a_{1,0}+a_{1,1}t+\ldots
+a_{1,n}t^{n},\ldots, a_{d,0}+a_{d,1}t+\ldots +a_{d,n}t^{n})$$
with coefficients in $k^{alg}$, such that
$f_{j}(a)=0\,mod\,t^{n+1}$ for each $j$. Jets can be considered as
approximate solutions for the system $f_j=0$. Using the $a_{k,l}$
as affine coordinates, we give the set of $n$-jets the structure
of a subvariety $\mathcal{L}_{n}(X)$ of $\A^{d(n+1)}_k$. There are
obvious truncation maps
$$\pi^{m}_{n}:\mathcal{L}_{m}(X)\rightarrow \mathcal{L}_{n}(X)$$
for $m\geq n$. Similarly, the set of arcs on $X$, that is,
$d$-tuples $\psi$ of power series over $k^{alg}$ satisfying
$f_j(\psi)=0$ for each $j$, can be seen as the set of closed
points of a $k$-scheme $\mathcal{L}(X)$, which comes with
truncation maps $$\pi_{n}:\mathcal{L}(X)\rightarrow
\mathcal{L}_{n}(X)\,.$$ Exact constructions and definitions are
given in the next section.

We can attach three motivic generating series to the variety $X$.
The Igusa Poincar\'e series $Q(T)$ counts all $n$-jets in
$\mathcal{L}_{n}(X)$, using the universal Euler characteristic,
taking values in the Grothendieck group of varieties over $k$. The
geometric Poincar\'e series $P_{geom}(T)$ only takes $n$-jets into
account which can be lifted through $\pi_{n}$ to an arc on $X$.
Finally, the arithmetic Poincar\'e series $P_{arith}$ counts, for
each field $K$ containing $k$, the $K$-rational $n$-jets that can
be lifted to a $K$-rational point of $\mathcal{L}(X)$. All three
series are rational (over the appropriate coefficient rings), see
the work of Denef and Loeser, in particular the survey article
\cite{DL2}.

While the Igusa Poincar\'e series can easily be expressed in terms
of a resolution of singularities, the geometric and arithmetic
series are very hard to compute in general. The proof of their
rationality is a qualitative proof, using results from model
theory like quantifier elimination, and does not yield
quantitative results. Up to now, the series had only been computed
for analytic branches of plane curves, and for toric surfaces
\cite{DL}\cite{LejReg}\cite{Nic1}, using Puiseux pairs and a
specific representation of arcs on toric varieties. In this paper,
we present a formula for both series in terms of a resolution of
singularities satisfying certain conditions, provided that such a
resolution exists. In fact, we will prove that in this case, the
series coincide. This opens a whole new realm of varieties $X$ for
which the series can easily be computed, including the toric
surfaces. In particular, our results imply the rationality of the
geometric and arithmetic series. Furthermore, our methods leave
much room for generalization, unlike the methods used to compute
the cases mentioned above. We believe that, at least in theory,
you can use similar arguments to compute the geometric Poincar\'e
series of any variety, the only restriction being the
combinatorial complexity.

To be precise, we determine the maximal truncation of an arc
$\psi$ in smooth ambient space that can be lifted to an arc on
$X$, and we construct an optimal approximation in
$\mathcal{L}(X)$, all in terms of the exact location of the
lifting of $\psi$ through the resolution morphism on the
exceptional locus.

Sections 2 and 3 contain some preliminaries on jets and motivic
integrals, and section 4 deals with the plane curve case. In
section 5, the general formulae for $P_{geom}$ are established.
Section 6 partially answers a question from \cite{LejReg},
concerning quasirational singularities. In section 7, we give a
very short computation of the geometric series of a toric surface.
Finally, in Section 8, we discuss the arithmetic series.

\section{Motivic integration and the geometric Poincar\'e series}
Until further notice, $k$ is an algebraically closed field of
characteristic zero.

Let $X$ be a variety over $k$, that is, a reduced and separated
scheme of finite type over $k$, not necessarily irreducible. For
each positive integer $n$, the functor from the category of
$k$-algebras to the category of sets, sending an algebra $R$ to
the set of $R[[t]]/t^{n+1}R[[t]]$-rational points on $X$, is
representable by a variety $\mathcal{L}_{n}(X)$. Since the natural
projections $\pi^{n+1}_{n}:\mathcal{L}_{n+1}(X)\rightarrow
\mathcal{L}_{n}(X)$ are affine, we can take the projective limit
in the category of schemes to obtain the scheme of arcs
$\mathcal{L}(X)$. This scheme represents the functor sending a
$k$-algebra $R$ to the set of $R[[t]]$-rational points on $X$, and
comes with natural projections $\pi_{n}:\mathcal{L}(X)\rightarrow
\mathcal{L}_{n}(X)$, mapping an arc to its $n$-truncation. For a
subvariety $Z$ of $X$, we define $\mathcal{L}(X)_Z$ to be the
closed subscheme $\pi_0^{-1}(Z)$ of $\mathcal{L}(X)$. When $X$ is
smooth, the morphisms $\pi_{n}^{n+1}$ are Zariski-locally trivial
fibrations with fiber $\A^{d}_{k}$, where $d$ is the dimension of
$X$. A morphism $h$ from $X$ to $Y$ induces a morphism $h$ from
$\mathcal{L}(X)$ to $\mathcal{L}(Y)$ by composition.

 We now introduce the Grothendieck ring $K_{0}(Var_k)$ of varieties over $k$.
 Start from the free abelian group generated by isomorphism
 classes $[X]$ of varieties $X$ over $k$, and consider the
 quotient by the relations $[X]=[X\setminus X']+[X']$, where $X'$
 is closed in $X$. A constructible subset of $X$ can be written as a disjoint union of
 locally closed subsets and determines unambiguously an element of $K_{0}(Var_k)$.
  The Cartesian product induces a product on
 $K_{0}(Var_k)$, which makes it a ring. We denote the class of the affine line
 $\A^{1}_{k}$ in $K_{0}(Var_k)$ by $\LL$, and the localization of $K_{0}(Var_k)$ with
 respect to $\LL$ by $\mathcal{M}_k$. On $\mathcal{M}_{k}$, we consider a decreasing filtration $F^{m}$, where $F^{m}$
 is the subgroup generated by elements of the form $[X]\LL^{-i}$, with dim\,$X-i\leq -m$. We define
 $\hat{\mathcal{M}}_{k}$ to be the completion of $\mathcal{M}_{k}$ with respect to this filtration.

 The
geometric Poincar\'e series, a formal power series over
$K_0(Var_k)$, is defined to be
$$P_{geom}(T)=\sum_{n\geq 0}[\pi_{n}(\mathcal{L}(X))]T^{n}\ .$$
Denef and Loeser \cite{DLinvent} proved that it is rational in
$\mathcal{M}_{k}[[T]]$.
 The
series is well defined, since Greenberg's theorem \cite{Gr} states
that we can find a positive integer $c$ such that, for all $n$,
and for each field $K$ containing $k$,
$\pi_{n}(\mathcal{L}(X)(K))=\pi_{n}^{nc}(\mathcal{L}_{nc}(X)(K))$.
So it follows from Chevalley's theorem \cite{Hart} that
$\pi_{n}(\mathcal{L}(X))$ is constructible, and hence determines
an element $[\pi_{n}(\mathcal{L}(X))]$ in $K_{0}(Var_k)$. One can
define local variants of this series, e.g. by only considering
arcs with origin in a fixed point $x$ of $X$.

Let $X\subset Y$ be varieties over $k$, with $Y$ smooth and of
dimension $d$. We define a valuation $d_{X}$ on $\mathcal{L}(Y)$
as follows: $d_{X}(\psi)=s$ if $\pi_{s-1}(\psi)\in
\pi_{s-1}\mathcal{L}(X)$, but $\pi_{s}(\psi)\notin
\pi_s\mathcal{L}(X)$, where we consider $\mathcal{L}(X)$ as a
subspace of $\mathcal{L}(Y)$. We define $d_{X}(\psi)$ to be
$\infty$ when $\psi\in \mathcal{L}(X)$, and to be $0$ if
$\psi(0)\notin X$. When $X$ is smooth, $\psi$ is a $k$-rational
arc on $Y$ with origin at $x$, and $\mathcal{I}$ is the defining
ideal sheaf of $X$ in $Y$,
$$d_{X}(\psi)=ord_t\,\mathcal{I}(\psi):=\min\{ord_t\,f(\psi)\,|\,f\in
\mathcal{I}_{x}\}\,.$$ In general, we will call
$ord_t\,\mathcal{I}(\psi)$ the order of contact between $\psi$ and
$X$, and we will denote this by $c(\psi,X)$.

For each positive integer $s$, and each point $x$ on $X$, we
define $D_{x}(X,s)$ by the motivic integral
$$D_x(X,s)=\int_{\mathcal{L}(Y)_x}\LL^{-d_{X}(\psi)s}d\mu(\psi)\,.$$
We refer the reader to
\cite{DLinvent}\cite{DL6}\cite{DL4}\cite{Nic1} for an introduction
to motivic integration. The normalization of the motivic measure
we use is the same as in these articles.

When it is clear which point $x$ and which variety $X$ the
integral $D_{x}(X,s)$ is associated to, we omit the subscript $x$
and the variable $X$ from our notation. Observe that $D_{x}(X,s)$
also depends on the ambient space $Y$.
 Putting $\LL^{-s}$ equal to $T$, the following simple formula
relates the local geometric Poincar\'e series $P_{geom}$ of $X$ at
$x$ to $D_{x}(s)$.
\begin{lemma}
$$P_{geom}(\LL^{-d}T)=\frac{1-\LL^{d}D_{x}(T)}{1-T}\  \mathrm{in}\  \hat{\mathcal{M}}_k[[T]]\,.$$
\end{lemma}
\begin{proof}
The general $T^{n}$-term of the right side can be written as
$$1-\LL^{d}\sum_{i=1}^{n}D_x(T)[i]\,$$ where $D_x(T)[i]$ denotes the
coefficient of $T^{i}$ in $D_x(T)$. Now it suffices to observe
that $P_{geom}(\LL^{-d}T)[n]$ equals $\LL^{d}$ times the motivic
measure of the arcs $\psi$ in $\mathcal{L}(X)_{x}$ satisfying
$d_{X}(\psi)>n$, while $1$ is equal to $\LL^{d}$ times the total
measure of $\mathcal{L}(X)_{x}$, and $\sum_{i=1}^{n}D_x(T)[i]$ is
the measure of the cylinder of arcs $\psi$ satisfying
$d_{X}(\psi)\leq n$. The lemma now follows from the additivity of
$\mu$.
\end{proof}
\section{Computing motivic integrals}\label{comp}
This section contains some trivial remarks, concerning the
computation of motivic integrals, using blow-ups and the change of
variables formula.

Let $X_1,\ldots,X_t$ be smooth subvarieties of a smooth
$d$-dimensional ambient space $Y$, intersecting transversally
along a subvariety $Z$. Let $c_i$ be the codimension of $X_i$ in
$Y$. We suppose that the sum of the $c_i$ does not exceed the
dimension $d$ of $Y$. One can check immediately that
\begin{equation}
\int_{\mathcal{L}(Y)_{Z}}\LL^{-\sum \alpha_i
c(\psi,X_i)}d\mu(\psi)=[Z]\LL^{-d}\prod_{i=1}^{t}(\LL^{c_i}-1)\frac{\LL^{-\alpha_i-c_i}}{1-\LL^{-\alpha_i-c_i}}\,,
\label{transv}
\end{equation}
where the coefficients $\alpha_i$ are positive integers. In
practice, one reduces to this situation using resolution of
singularities and the change of variables formula for motivic
integrals.

Now suppose that $t=2$, and let $A$ be the measurable subset of
$\mathcal{L}(Y)$ consisting of all arcs $\psi$, satisfying
$\psi(0)\in Z$ and $c(\psi,X_1)\geq c(\psi,X_2)$. We can easily
compute the motivic integral
$$I=\int_{A}\LL^{-\alpha_1c(\psi,X_1)-\alpha_2c(\psi,X_2)}d\mu(\psi)$$
by blowing up $Z$, and applying the change of variables formula.
Let $\psi$ be an arc on $Y$, with origin in $Z$ but not entirely
contained in $Z$, and let $\psi'$ be its lifting through the
blow-up $h:Y'\rightarrow Y$ with center $Z$ and exceptional
divisor $E$. Let $X'_i$ be the strict transform of $X_i$. It is
clear that $\psi$ belongs to $A$ if and only if $\psi'(0)\notin
X'_2$. The change of variables formula yields
\begin{eqnarray*}
I&=&\int_{\mathcal{L}(Y')_{E\setminus
X'_2}}\LL^{-\alpha_1c(\psi',X'_1)-(\alpha_1+\alpha_2+c_1+c_2-1)c(\psi',E)}d\mu(\psi')
\\&=&\LL^{-d}(\LL-1)\{[E\setminus (X'_1\cup
X'_2)]\frac{\LL^{-(\alpha_1+\alpha_2+c_1+c_2)}}{1-\LL^{-(\alpha_1+\alpha_2+c_1+c_2)}}
\\[5pt]&&+[E\cap X'_1](\LL^{c_1}-1)\frac{\LL^{-(2\alpha_1+\alpha_2+2c_1+c_2)}}{(1-\LL^{-(\alpha_1+\alpha_2+c_1+c_2)})(1-\LL^{-\alpha_1-c_1})}\}\,.
\end{eqnarray*}

Of course, analogous statements can be formulated for $t>2$, or
more complicated sets $A$.

As a final example, we compute the motivic integral
$$I=\int_{\mathcal{L}(Y)_{Z}}\LL^{-\lfloor
c(\psi,X_1)/2+c(\psi,X_2)/2\rfloor }d\mu(\psi)\,,$$ where $\lfloor
x\rfloor$ is the largest integer smaller than or equal to $x$.
Blowing up $Z$, as before, we get
$$I=\int_{\mathcal{L}(Y')_{E}}\LL^{-(c_1+c_2)c(\psi',E)-\lfloor
c(\psi',X'_1)/2+c(\psi',X'_2)/2\rfloor}d\mu(\psi')\,.$$ The
advantage of this method is that $X'_1\cap E$ and $X'_2\cap E$ are
disjoint. A straightforward computation yields
\begin{eqnarray*}
\lefteqn{I=\LL^{-d}(\LL-1)\frac{\LL^{-(c_1+c_2+1)}}{1-\LL^{-(c_1+c_2+1)}}\{[E\setminus
(X'_1\cup X'_2)]}
\\&&\qquad+[E\cap
X'_1](\LL^{c_1}-1)\frac{\LL^{-c_1}+\LL^{-2c_1-1}}{1-\LL^{-2c_1-1}}
+[E\cap
X'_2](\LL^{c_2}-1)\frac{\LL^{-c_2}+\LL^{-2c_2-1}}{1-\LL^{-2c_2-1}}\}\,.
\end{eqnarray*}
 These and similar methods will allow us to compute the
geometric Poincar\'e series from the formula for $d_X$, that we
will establish in a subsequent section. As horrifying as the
computations may look, they are all based on the same basic
principles: blowing up in order to obtain transversal
intersection, and to simplify the integration domain, and applying
the change of variables formula.
\section{Plane curves}
 Let $k$ be an algebraically closed field of
characteristic 0. Let $X$ be a formal branch of a plane algebraic
curve over $k$, with Puiseux expansion
$$\left\{\begin{array}{lll}
x&=&t^{m}
\\y&=&\sum_p a_pt^{p}\,,
\end{array}\right.$$
where we suppose $m<\min\{p\,|\,a_p\neq 0\}$. Let
$(p_i,q_i),\,i=1,\ldots,s$ be the characteristic pairs, and let
$c_1t^{k_1}+\ldots+c_st^{k_s}$ be the corresponding essential
terms in the expansion of $y$. As is explained in \cite{Bries},
the knowledge of the characteristic pairs suffices to construct
the resolution graph of the minimal embedded resolution
$h:\tilde{Y}\rightarrow Y=\A^{2}_k$ of $X$. We will derive a
formula for the valuation function $d_X\circ h$ in terms of the
contact of an arc on $\tilde{Y}$ with the exceptional locus of
$h$, and with the strict transform $\tilde{X}$ of $X$.

Let us recall what the resolution graph of $h$ looks like. In
addition to the numerical data $(N_j,\nu_j)$, we attach to each
irreducible exceptional divisor $E_j$ a couple $(M_j,\mu_j)$. Let
$E_k$, $k\in K(j)$, be the exceptional divisors containing the
point that was blown up in the creation of $E_j$; thus $K(j)$ has
at most two elements, and may be empty. We define the multiplicity
$M_j$ to be equal to $N_j$, subtracted by the sum of the $N_k$,
$k\in K(j)$, and we put $\mu_j$ equal to the sum of the $\mu_k$,
where the divisor emerging in the first blow-up gets initial value
$\mu=1$.

 For each index $i=1,...,s$, we get a
chain $P_i$ of exceptional divisors. In order to describe what
$P_i$ looks like, we first introduce some new invariants. Consider
the following Euclidean algorithm:
\begin{eqnarray*}
\kappa_i&=&a_{i,1}r_{i,1}+r_{i,2}
\\r_{i,1}&=&a_{i,2}r_{i,2}+r_{i,3}
\\&&\ldots
\\r_{i,w(i)-1}&=&a_{i,w(i)}r_{i,w(i)}\,,
\end{eqnarray*}
where $\kappa_1=k_1$, $\kappa_i=k_i-k_{i-1}$, $r_{1,1}=m$, and
$r_{i,1}=r_{i-1,w(i-1)}$. Given $P_1,\ldots,P_{i-1}$, the
resolution process runs as follows: first, we get a chain of
divisors $E_{i,1,j}$, $j=1,\ldots,a_{i,1}$, each with multiplicity
$r_{i,1}$. Thereupon, $a_{i,2}$ divisors $E_{i,2,j}$ with
multiplicity $r_{i,2}$ emerge, each of them separating the
previous one from $E_{i,1,a_{i,1}}$. This process continues, so at
the end, $P_i$ is the chain starting with
$$E_{i,1,1},\ldots,E_{i,1,a_{i,1}},E_{i,3,1},\ldots,E_{i,3,a_{i,3}},\ldots$$
and ending in
$$\ldots,E_{i,4,a_{i,4}},\ldots,
E_{i,4,1},E_{i,2,a_{i,2}},\ldots,E_{i,2,1}\ .$$ Denote
$E_{i,w(i),a_{i,w(i)}}$ by $F_i$. If $i<s$, the divisor
$E_{i+1,1,1}$ intersects $F_{i}$ in a smooth point of $P_{i}$, and
$\tilde{X}$ intersects $F_{s}$ in a smooth point of $\cup P_{i}$.

Let $\psi$ be a non-constant $k$-rational arc on $\A^{2}_{k}$,
with $\psi(0)=(0,0)$, and let $\tilde{\psi}$ be the unique lifting
of $\psi$ through $h$. By a lifting $\varphi'$ of an arc $\varphi$
on a variety $Z$ through a proper birational morphism
$f:W\rightarrow Z$, we mean the following: suppose that $f$ is an
isomorphism over a Zariski-open neighbourhood in $Z$ of the
generic point of $\varphi$. Applying the valuative criterion for
properness to the morphism $f$, we see that there exists a unique
arc $\varphi'$ on $W$ such that the composition $h\circ \varphi'$
is equal to $\varphi$. We call this $\varphi'$ the lifting of
$\varphi$ through $f$.

 We will give a formula for $d_{X}(\psi)$ in
terms of the contact of $\tilde{\psi}$ with the exceptional locus
and $\tilde{X}$. This is possible because this contact information
allows us to reconstruct the relevant part of the power series
expansion of $\psi$.

 If $x=\tilde{\psi}(0)$ is a smooth point of $\cup P_{i}$, contained in $E_{1,1,1}$, the constant arc at the origin of
$\A^{2}_{k}$ is an optimal approximation for $\psi$ in
$\mathcal{L}(X)$ (with respect to the valuation $ord_t$), so
$d_X(\psi)$ equals the multiplicity $n$ of $\psi$. If not,
 we may assume $\psi(s)$ to be of the form
$$\left\{\begin{array}{lll}
x&=&s^{n}
\\y&=&\sum_p b_ps^{p}\,,
\end{array}\right.$$
with $n<\min\{p\,|b_p\neq 0\}$.
 In both cases, it is easily
verified that the multiplicity $n$ is equal to $\sum
\mu_{k}\gamma_k$, where the sum is taken over the exceptional
components $E_k$ containing $x$, and $\gamma_k$ is the order of
contact of $\tilde{\psi}$ with $E_k$. If $n$ is not a multiple of
$m$, the constant arc at the origin will again be an optimal
approximation for $\psi$, so $d_X(\psi)=n$. From now on, we
suppose that $n=\lambda m$, $\lambda\in \N$.

First, suppose that $x\notin \tilde{X}$ is contained in
$E_{i,j,k}\neq F_i$, with $j$ even. This implies that $\psi$
agrees with the Puiseux expansion of $X$, modulo a
reparametrization $t=s^{\lambda}$, up to the essential term
$c_it^{k_i}$, so $d_{X}(\psi)=\lambda k_i$. If $j$ is odd,
$E_{i,j,k}$ is the only exceptional component containing $x$, and
$x\notin \tilde{X}$, it is easy to see that
$$\lambda^{-1}d_{X}(\psi)=k_{i-1}+a_{i,1}r_{i,1}+a_{i,3}r_{i,3}+\ldots+k\,r_{i,j},$$
 where we put $k_0=0$. Next, assume that $x\notin \tilde{X}$ is
the intersection point of $E_{i,j,k}$ and $E_{i,j,k+1}$, where $j$
is odd. In this case,
$$\lambda^{-1}d_{X}(\psi)=k_{i-1}+a_{i,1}r_{i,1}+a_{i,3}r_{i,3}+\ldots+k\,r_{i,j}+\lambda^{-1}\gamma,$$
where $\gamma$ is the order of contact of $\tilde{\psi}$ with
$E_{i,j,k+1}$. If $x$ is the intersection point of
$E_{i,j,a_{i,j}}$ and $E_{i,j+2,1}$, where $j$ is odd, then
$$\lambda^{-1}d_{X}(\psi)=k_{i-1}+a_{i,1}r_{i,1}+a_{i,3}r_{i,3}+\ldots+a_{i,j}\,r_{i,j}+\lambda^{-1}\gamma,$$
where $\gamma$ is the order of contact of $\tilde{\psi}$ with
$E_{i,j+2,1}$. And if $x$ is the intersection point of
$E_{i,j,a_{i,j}}$ and $F_i$, where $j$ is odd, and $w(i)$ even,
$$\lambda^{-1}d_{X}(\psi)=k_{i-1}+a_{i,1}r_{i,1}+a_{i,3}r_{i,3}+\ldots+a_{i,j}\,r_{i,j}+\lambda^{-1}\gamma,$$
where $\gamma$ is the order of contact of $\tilde{\psi}$ with
$F_i$.
 Finally, if $x\in \tilde{X}$, $d_{X}(\psi)$
equals $\lambda k_{s}+\gamma$, where this time $\gamma$ is the
contact order of $\tilde{\psi}$ with $\tilde{X}$.

 This analysis
allows one to compute the geometric Poincar\'e series of $Y$ -
which was already computed in a much more elementary way in
\cite{DL} - using the motivic change of variables formula
$$\int_{\mathcal{L}(Y)}\LL^{-d_{X}}d\mu=\int_{\mathcal{L}(\tilde{Y})}\LL^{-d_{X}\circ
h - ord_t\, Jac_{h}}d\mu\,.$$ One can use the same techniques to
compute the geometric Poincar\'e series for plane curves which are
not necessarily analytically irreducible.

\section{General results}
Let $X\subset Y$ be varieties over an algebraically closed field
$k$ of characteristic zero, with $Y$ smooth and of dimension $d$,
and $X$ of dimension $m$. Let $h:Y'\rightarrow Y$ be a composition
$h_r\circ\dots\circ h_1$, where each $h_i$ is the blow-up of a
point, and, if $i>1$, this point lies on at most one exceptional
divisor of $h_{i-1}\circ\dots\circ h_1$. Let $X'$ be the strict
transform of $X$ under $h$, and suppose that $X'$ is smooth.
Assume furthermore that $X'$ and the exceptional locus $E$
intersect transversally, and that each exceptional component $E_i$
of $E$ contains a point of $X'$ that does not lie on any other
exceptional component. Let $\psi$ be $k$-rational arc in
$\mathcal{L}(Y)_{X}\setminus \mathcal{L}(X)$, and let $\psi'$ be
the lifting of
$\psi$ through $h$. 
 If $\psi'(0)\notin E$, put $\lambda=0$. If $\psi'(0)$ lies on
precisely one exceptional component $E_i$ of $E$, we define
$\lambda$ to be equal to $c(\psi',E_i)$, multiplied by the order
$\nu_i$ of the Jacobian of $h$ on $E_i$, divided by $d-1$. This
latter factor $\nu_i/(d-1)$ indicates the depth $e_i$ of $E_i$ in
the composition of blow-ups. Finally, if $\psi'(0)$ lies on two
exceptional components $E_{i+1}$ and $E_i$, where $E_{i+1}$ was
created by blowing up a point of $E_i$, we put $\lambda$ equal to
$$c(\psi',E_i)e_i+c(\psi',E_{i+1})e_{i+1}\,.$$ With
this notation, we can formulate the following theorem.
\begin{theorem}\label{formula}
Under the conditions explained above,
$$d_{X}(\psi)=d_{X'}(\psi')+\lambda\,.$$
\end{theorem}
\begin{proof}
The fact that we are considering arcs, allows us to work locally
with respect to the \'etale topology.
 Using our assumptions on $h$, the following lemma is
easily verified.

\begin{lemma}\label{formh}
We can find local coordinates $(y_1,\ldots,y_{d})$ on $Y$ at
$\psi(0)$, and local coordinates $(y'_1,\ldots,y'_{d})$ on $Y'$ at
$\psi'(0)$, such that the following properties are satisfied:
\begin{itemize}
\item If $\psi'(0)$ is contained in exactly one exceptional
component $E_i$, the morphism $h$ is given by
$$h(y'_1,\ldots,y'_d)=(y'_1,(y'_1)^{e_i}y'_2,\ldots,(y'_1)^{e_i}y'_m, (y'_1)^{e_i}\Phi_{m+1},\ldots ,(y'_1)^{e_i}\Phi_d)\,.$$
 If $\psi'(0)$ is
contained in two distinct components $E_i$ and $E_{i+1}$,
$h(y'_1,\ldots,y'_d)$ is equal to
$$\!\!(y'_1y'_2,(y'_1)^{e_i}(y'_2)^{e_i+1},(y'_1)^{e_i}(y'_2)^{e_i+1}y'_3,
\ldots,(y'_1)^{e_i}(y'_2)^{e_i+1}y'_m,$$
\\[-25pt]$$\qquad\qquad\qquad\qquad\qquad\qquad\qquad (y'_1)^{e_i}(y'_2)^{e_i+1}\Phi_{m+1},\ldots,(y'_1)^{e_i}(y'_2)^{e_i+1}\Phi_d)\,.$$
Here $\Phi=(\Phi_{m+1},\ldots,\Phi_d)$ is a tuple of power series
in the maximal ideal $M$ of $k[[y'_1,\ldots,y'_d]]$, such that the
determinant of $[\frac{\partial \Phi_i}{\partial
y'_j}]_{i,j=m+1}^{d}$ is a unit at the origin. Furthermore, in the
first case, each $\Phi_j$ is contained in the ideal
$(y'_{j})+M^{2}$, and in the second case, each $\Phi_j$ is
contained in $(y'_1,y'_2,y'_{j})+M^{2}$.
 \item Either $\psi'(0)\notin X'$, and in
this case we can choose $\Phi_j=y'_j$ for all $j$; or $X'$ is
locally defined by
$$y'_{m+1}=\ldots=y'_d=0\,.$$
 \end{itemize}
 \end{lemma}
 \begin{proof}[Proof of the lemma]
We prove the lemma when $\psi'(0)\in E_i\cap E_{i+1}$; the other
case is easier. Choose local coordinates $(y_i)$ on $Y$ at
$\psi(0)$ such that the order of $y_1(\psi)$ is minimal among
$\{ord_t y_j(\psi)\}$. It follows from Hensel's lemma that we can
reparametrize $\psi$, that is, compose $\psi$ with an automorphism
of Spec\,$k[[t]]$, such that $y_1(\psi)=t^{c}$ for some positive
integer $c$. After choosing new coordinates on our ground space
$Y$, we may assume that $c$ does not divide the order of
$y_j(\psi)$ if $j\neq 1$. Let $c'$ be the minimum of $\{ord_t
y_j(\psi)\,|\,j=2,\ldots,d\}$. We may assume, changing coordinates
if necessary, that $ord_t y_j(\psi)=c'$ iff $j=2$. Since we never
blow up points that belong to two distinct exceptional components,
it is clear that $e_i$ is equal to $\lfloor c/c'\rfloor$, and that
we can find local coordinates $(z'_i)$ at $\psi'(0)$ such that $h$
is given by
$$h(z'_1,\ldots,z'_d)=(z'_1z'_2,(z'_1)^{e_i}(z'_2)^{e_i+1},(z'_1)^{e_i}(z'_2)^{e_i+1}z'_3,
\ldots,(z'_1)^{e_i}(z'_2)^{e_i+1}z'_d)\,.$$ Since $X',\,E_i$ and
$E_{i+1}$ intersect transversally at $\psi'(0)$, we can choose,
after permutating the $z'_j$ if necessary, new local coordinates
$(y'_i)$ at $\psi'(0)$, with $y'_j=z'_j$ for $j=1,\ldots,m$, such
that $X'$ is defined by $y'_{m+1}=\ldots=y'_d=0$. The part about
the $\Phi_j$ is obvious.
 \end{proof}

 If $\psi'(0)\in X'$, it is clear that $\varphi=h(\varphi')$ is an optimal approximation for $\psi$ in
 $\mathcal{L}(X)$, where $y'_i(\varphi')=y'_i(\psi')$ for
 $i=1,\ldots,r-1$, and $y'_i(\varphi')=0$ if $i\geq r$, and that
 the postulated formula for $d_X(\psi)$ holds. For suppose that
 $\eta$ is an optimal approximation, and that $\eta$ is a better approximation than $\varphi$. The fact that $\eta$
 lies at least as close to $\psi$ as $\varphi$ does (with respect to the valuation $ord_t$), guarantees that the
 lifting $\eta'$ of $\eta$ through $h$ has its origin at $\psi'(0)$. It is
 clear that $y'_j(\eta')=0$ for $j=m+1,\ldots,d$. Suppose that the
 minimum of $\{ord_t y'_j(\psi')\,|\,j=m+1,\ldots,d\}$ is
 realized for $j=d$.

 If $\psi'(0)$ is contained in exactly one
 exceptional component, the fact that $\eta$ is a better
 approximation than $\varphi$ implies that
 $$y'_1(\psi')\equiv y'_1(\eta')\,mod\, t^{e_i ord_t y'_1(\psi')+ord_t y'_d(\psi')+1}\,,$$
  and that
$$\Phi_d(y'_1(\psi'),\ldots,y'_d(\psi'))\equiv \Phi_d(y'_1(\eta'),\ldots,y'_d(\eta'))\,mod\,t^{ord_t y'_d(\psi')+1}\,.$$
 Since $\Phi_d\in (y'_d)+M^{2}$, and $\frac{\partial
 \Phi_d}{\partial y'_d}\neq 0$, and $y'_d(\eta')=0$, we see that
$ ord_t (y'_j(\eta')-y'_j(\psi'))<ord_t y'_d(\psi')$ for some
$j\in
 \{2,\ldots,m\}$, which makes $\varphi$ a better approximation
 than $\eta$.

If $\psi'(0)\in E_i\cap E_{i+1}$, we can follow similar arguments:
we see that
\begin{eqnarray*}
 y'_1(\psi')y'_2(\psi')&\equiv&
y'_1(\eta')y'_2(\eta')\,mod\, t^{e_i ord_t
y'_1(\psi')+e_{i+1}ord_t y'_2(\psi')+ ord_t y'_d(\psi')+1}\,,
\\y'_1(\psi')^{e_i}y'_2(\psi')^{e_i+1}&\equiv& y'_1(\eta')^{e_i}y'_2(\eta')^{e_i+1}\,mod\, t^{e_i
ord_t y'_1(\psi')+e_{i+1}ord_t y'_2(\psi')+ ord_t
y'_d(\psi')+1}\,,
\\ \Phi_d(y'_1(\psi'),\ldots,y'_d(\psi'))&\equiv&
\Phi_d(y'_1(\eta'),\ldots,y'_d(\eta'))\,mod\,t^{ord_t
y'_d(\psi')+1}\,.
\end{eqnarray*}
These observations again lead to the conclusion that $ ord_t
(y'_j(\eta')-y'_j(\psi'))<ord_t y'_d(\psi')$ for some $j\in
 \{2,\ldots,m\}$ (use the fact that for units $u,v$ in $k[[t]]$,
 the congruence $u\equiv v \,mod\,t^{a}$ implies $u^{-1}\equiv v^{-1}\,mod\,t^{a}$, where $a\in \N$) .

 Now assume that $\psi'(0)\notin X'$, and that $\psi'(0)$ is
 contained in exactly one exceptional component $E_i$. Let $h_i$ be the blow-up of the point $x$, creating $E_i$, and
 decompose $h$ as $h=\tilde{h}\circ h_i\circ h'$, where $h$ and $h'$ are compositions of blow-ups.
It follows from the structure of the resolution $h$ of $X$ that
there exists an arc $\varphi'$ on $X'$ such that $h_i\circ
h'\circ\varphi'(0)=x$ and $y_1\circ h(\varphi')=y_1(\psi)$, and
such that $E_i$ is the only exceptional component containing
$h'\circ \varphi'(0)$. Let $\eta$ be a non-constant optimal
approximation of $\psi$ in $\mathcal{L}(X)$. An argument similar
to
 the one used above, shows that the lifting of $\eta$ through $\tilde{h}$
 has its origin at $x$. If the distance from $\eta$ to $\psi$ were
 strictly smaller than the distance from $\varphi=h(\varphi')$ to
 $\psi$, this would imply that the lifting of $\eta$ through
 $\tilde{h}\circ h_i$ would lie in $X'$, which contradicts our
 assumptions. Hence, the formula for $d_{X}(\psi)$ holds in
 this case also. An analogous reasoning can be used when $\psi'(0)$ is contained in $E_i\cap E_{i+1}\setminus X$.
\end{proof}

%

Using this formula, and the expression \ref{transv} in Section
\ref{comp}, it becomes very easy to compute the geometric
Poincar\'e series for singularities with an embedded resolution
satisfying the conditions of Theorem \ref{formula} .

Now suppose that $X$ is a surface, and all conditions for Theorem
\ref{formula} are satisfied, except that we allow $E$ and $X'$ to
intersect non-transversally in a finite number of points $x_i$,
each of which has to be a smooth point of $E$, at which the
intersection multiplicity is two. We suppose that locally around
$x_i$ (with respect to the \'etale topology), the intersection of
$E$ and $X'$ consists of two smooth prime divisors $F_1$ and
$F_2$, meeting transversally at $x_i$ (the obvious generalizations
hold when the pair $(X',E)$ is locally a product of a pair of the
required form with a smooth space). The proof of Theorem
\ref{formula} remains valid for each arc $\psi$ whose lifting
$\psi'$ does not have its origin at one of the $x_i$.

 We obtain
transversal intersection by blowing up one of them, say $F_1$. Let
$h': \tilde{Y}\rightarrow Y'$ be the blow-up morphism, with
exceptional divisor $E'_i$, and let $\tilde{X}$ be the strict
transform of $X'$. Suppose that the arc $\psi$ lifts to an arc
$\tilde{\psi}$, with origin on $E$, through the composition
$h'\circ h$. Put $\psi'$ equal to $h'(\tilde{\psi})$. Let
$\tilde{F}$ be the inverse image of $x_i$ under the projection
$E'_i\cong F_1\times \Pro_{k}^{d-2}\rightarrow F_1$. Blowing up
$F_2$ yields an exceptional divisor $E''_i$, and a strict
transform $\bar{X}$. Denote by $\bar{\psi}$ the lifting of
$\tilde{\psi}$ through this
 blow-up morphism.

Let $\tilde{x}_i$ be the intersection of $\tilde{X}$ with
$\tilde{F}$. By arguments similar to the ones in Lemma
\ref{formh}, we can find local coordinates
$(\tilde{y}_1,\ldots,\tilde{y}_d)$ on $\tilde{Y}$ at
$\tilde{x}_i$, local coordinates $(y'_1,\ldots,y'_d)$ on $Y'$ at
$x_i$, and local coordinates $(y_1,\ldots,y_d)$ on $Y$ at
$\psi(0)$, such that
\begin{itemize}
\item the morphism $h$ is given in local coordinates by
$$h(y'_1,\ldots,y'_d)=(y'_1,(y'_1)^{e_i}\Phi_2,\ldots,(y'_1)^{e_i}\Phi_d)\,,$$
where $\Phi_j\in (y_j)+M^{2}$ for each $j$, with $M$ as in Lemma
\ref{formh}.
 \item the morphism $h'$ is given in local coordinates
by
$$h'(\tilde{y}_1,\ldots,\tilde{y}_d)=(\tilde{y}_1\tilde{y}_2,\tilde{y}_2,\tilde{y}_3,\tilde{y}_2\tilde{y}_4,\ldots,\tilde{y_2}\tilde{y}_d)\,,$$
\item the strict transform of $E_i$ under $h'$ is defined by
$\tilde{y}_1=0$, \item $\tilde{X}$ is defined by
$\tilde{y}_3=\tilde{y}_1$ and $\tilde{y}_4=\ldots=\tilde{y}_d=0$.
\item The smooth germ $X'$ is locally defined by the equations
$y'_1-y'_2\,y'_3=0$ and $y'_4=\ldots =y'_d=0$.
\end{itemize}
Let $H'_i$ be locally defined by $y'_2=0$, and $H''_i$ by
$y'_3=0$. We use the same notation for their strict transforms
under any blow-up.

First, we suppose $d=3$. It follows from the expression for $h$,
and an argument similar to the one used in the proof of Theorem
\ref{formula}, that we find an optimal approximation
$\varphi=h(\varphi')$ by maximizing
\begin{eqnarray*}
\lefteqn{d(\psi',\varphi'):=\min\,\{ord_t\,(y'_1(\varphi')-y'_1(\psi')),\,
ord_t\,((y'_1)^{e_i}(\varphi')y'_2(\varphi')-(y'_1)^{e_i}(\psi')y'_2(\psi')),\,}\qquad\qquad\qquad\qquad\qquad\qquad\qquad\\
&&ord_t\,((y'_1)^{e_i}(\varphi')y'_3(\varphi')-(y'_1)^{e_i}(\psi')y'_3(\psi'))\}\,.\end{eqnarray*}

After a suitable reparametrization, we may suppose that
$y'_1(\psi')=t^{N}$. If the leading terms of $y'_1(\psi')$ and
$y'_1(\varphi')$ differ, it is clear that we can replace
$\varphi'$ by another arc on $X'$ whose image lies at least as
close to $\psi$ as $\varphi$ does, whose leading $y'_1$-term
agrees with that of $\psi'$. Hence, we might as well assume that
they agreed all along.
\begin{lemma}
 It suffices to maximize
\begin{eqnarray*}
\lefteqn{\min\,\{ord_t\,(y'_1(\varphi')-y'_1(\psi')),\,
e_i\,ord_t\,y'_1(\psi')+ ord_t\,(y'_2(\varphi')-y'_2(\psi')),\,}\qquad\qquad\qquad\qquad\qquad\qquad\qquad\\
&&e_i\,ord_t\,y'_1(\psi')+
ord_t\,(y'_3(\varphi')-y'_3(\psi'))\}\,.\end{eqnarray*} This can
always be achieved without modifying $y'_1$, i.e. we can put
$y'_1(\varphi')=y'_1(\psi')$.
\end{lemma}
\begin{proof}[Proof of the lemma]
Given an approximation $\varphi'$ in $\mathcal{L}(X')_{x_i}$, with
$y'_1(\varphi')\neq y'_1(\psi')=t^{N}$, we construct an arc
$\eta'$ in $\mathcal{L}(X')_{x_i}$, such that $y'_1(\eta')=t^{N}$,
and such that $d(\eta',\psi')\geq d(\varphi',\psi')$.

We already noticed that we may suppose
$y'_1(\varphi')=t^{N}+t^{N+1}\phi'$. Now define an arc $\eta'$ by
$y'_1(\eta')=t^{N}$, $y'_2(\eta')=y'_2(\varphi')(1+t\phi')^{-1}$,
and $y'_3(\eta')=y'_3(\varphi')$. Since $d(\psi',\varphi')$ is at
most $N+1+ord_t\,\phi'$, it is sufficient to show that
$d(\eta',\varphi')$ is at least $N+1+ord_t\,\phi'$. However, this
is clear from the definition.
\end{proof}

If $\tilde{\psi}(0)\notin \tilde{X}$, the order of $y'_2(\psi')$
is larger than, or equal to the order of $y'_1(\psi')$.
 If $ord_t\,y'_3(\psi')$ is smaller than
$\lceil ord_t\,y'_1(\psi')/2 \rceil$, we define $y'_2(\varphi')$
to be $y'_1(\psi')(y'_3(\psi'))^{-1}$. In the other case, we take
for $y'_2(\varphi')$ and $y'_3(\varphi')$ two power series of
order $\lfloor ord_t\,y'_1(\psi')/2 \rfloor$, resp. $\lceil
ord_t\,y'_1(\psi') \rceil/2$, whose product equals $y'_1(\psi')$.
So
$$d_{X}(\psi)= 
e_ic(\tilde{\psi},E'_i)+\max\,\{\lfloor 
c(\tilde{\psi},E'_i)/2\rfloor,\,
c(\tilde{\psi},E'_i)-c(\tilde{\psi},H''_i)\} \,.$$

Suppose that $\tilde{\psi}(0)$ is contained in $\tilde{X}$. If the
leading terms of $y'_1(\psi')$ and $y'_2(\psi')y'_3(\psi')$
differ, analogous arguments yield
\begin{eqnarray*}
\lefteqn{d_X(\psi)=e_ic(\psi',E_i)+\max\,\{\min\{\lfloor
c(\psi',E_i)/2\rfloor,c(\psi',H'_i) ,c(\psi',H''_i)\},}\quad
\\&& \min\{c(\psi',H'_i), c(\psi',E_i)-c(\psi',H''_i)\},\,\min\{c(\psi',H''_i),\,c(\psi',E_i)-c(\psi',H'_i)\}\,\}\,.
\end{eqnarray*}

 When the leading terms of
$y'_1(\psi')$ and $y'_2(\psi')y'_3(\psi')$ coincide,
\begin{eqnarray*}
d_{X}(\psi)&=&e_iord_t\,y'_1+d_{X'}(\psi')-\min\,\{c(\psi',H'_i),\,c(\psi',H''_i)\}
\\&=&e_ic(\tilde{\psi},E_i)+(e_i+1)c(\tilde{\psi},E'_i)+d_{\tilde{X}}(\tilde{\psi})-c(\tilde{\psi},\tilde{F})\,.
\end{eqnarray*}
Blowing up $F_2$, we see that the following formula holds in
general:
\begin{eqnarray*}\label{dix} \lefteqn{d_{X}(\psi)=d_{\bar{X}}(\bar{\psi})+ e_ic(\bar{\psi},E_i)+
e_ic(\bar{\psi},E'_i)+e_ic(\bar{\psi},E''_i)+\qquad\qquad\qquad\quad
}\\&&\nonumber+\max\,\{\,\min\{\lfloor c(\bar{\psi},E_i)/2+
c(\bar{\psi},E'_i)/2+c(\bar{\psi},E''_i)/2\rfloor,
\\&&\qquad\qquad\qquad\qquad\qquad\qquad c(\bar{\psi},E'_i)+c(\bar{\psi},H'_i),c(\bar{\psi},E''_i)+c(\bar{\psi},H''_i)\},
\\&&\min\{c(\bar{\psi},E'_i)+c(\bar{\psi},H'_i),\,c(\bar{\psi},E_i)+c(\bar{\psi},E'_i)-c(\bar{\psi},H''_i)\},
\\&&\min\{c(\bar{\psi},E''_i)+c(\bar{\psi},H''_i),\,c(\bar{\psi},E_i)+c(\bar{\psi},E''_i)-c(\bar{\psi},H'_i)\}\,\}\,.\qquad\qquad\quad (2)\end{eqnarray*}

Now suppose $d>3$. Let $\bar{Z}$ be the strict transform of the
\'etale germ $Z'$ at $x_i$, locally defined by $y'_1-y'_2y'_3=0$,
and let $W'$ be the germ at $x_i$ defined by $y'_4=\ldots
=y'_d=0$. Provided that we replace $d_{\bar{X}}(\bar{\psi})$ by
$d_{\bar{Z}}(\bar{\psi})$, our formula remains valid, unless
$e_ic(E_i)+ord_t\,y'_j$ is smaller than the above expression
(\ref{dix}), for some $j$ in $\{4,\ldots,m\}$. If this is the
case, $d_{X}(\psi)$ will be equal to
$e_ic(\psi',E_i)+d_{W'}(\psi')$.

We've proven
\begin{theorem}\label{form2}
The valuation $d_{X}(\psi)$ is the minimum of
\begin{eqnarray*}\label{nontransv}
\lefteqn{d_{\bar{Z}}(\bar{\psi})+ e_ic(\bar{\psi},E_i)+
e_ic(\bar{\psi},E'_i)+e_ic(\bar{\psi},E''_i)+\qquad\qquad\qquad\quad
}\\&&\nonumber+\max\,\{\,\min\{\lfloor c(\bar{\psi},E_i)/2+
c(\bar{\psi},E'_i)/2+c(\bar{\psi},E''_i)/2\rfloor,
\\&&\qquad\qquad\qquad\qquad\qquad\qquad c(\bar{\psi},E'_i)+c(\bar{\psi},H'_i),c(\bar{\psi},E''_i)+c(\bar{\psi},H''_i)\},
\\&&\min\{c(\bar{\psi},E'_i)+c(\bar{\psi},H'_i),\,c(\bar{\psi},E_i)+c(\bar{\psi},E'_i)-c(\bar{\psi},H''_i)\},
\\&&\min\{c(\bar{\psi},E''_i)+c(\bar{\psi},H''_i),\,c(\bar{\psi},E_i)+c(\bar{\psi},E''_i)-c(\bar{\psi},H'_i)\}\,\}\,.
\end{eqnarray*} and $d_{W'}(\psi')+e_ic(\psi',E_i)$.
\end{theorem}

 Let us compute the motivic integral
\begin{eqnarray*}
D^{(e_i)}(X,s)& = & \int_{\mathcal{L}(Y')_{x_i}}\LL^{-d_{X}\circ
h(\psi')s-ord_t\,Jac_h}d\mu(\psi')
\\ & = & \int_{\mathcal{L}(\tilde{Y})_{h'^{-1}(x_i)}}\LL^{-d_{X}\circ
h\circ h'(\tilde{\psi})s-ord_t\,Jac_{h\circ
h'}}d\mu(\tilde{\psi})\,.
\end{eqnarray*}

 You can visualize the situation
by taking $E_i'$ and $E''_i$ to be two walls making a rectangular
corner, and by imagining $\bar{X}$, $E_i$ and $H''_i$ to be
horizontal shelves. The fiber over $\tilde{x}_i$ is equal to
$(E'_i\cap E''_i)$. To simplify notation, we denote $e_i(d-1)+1$
by $\nu_1$, and $(e_i+1)(d-1)$ by $\nu_2$. The order of the
Jacobian of the resolution morphism on $E_i$ is $\nu_1$, and on
$E'_i$, as well as on $E''_i$, it equals $\nu_2$.

First, suppose $d=3$. The easiest way to compute $D^{(e_i)}$ would
be to take a shortcut and use the formula in \cite{LejReg} for the
geometric Poincar\'e series for toric surfaces. The minimal
resolution of the toric hypersurface defined by $xz-y^{2e_i+1}=0$
has a point of non-transversal intersection at depth $e_i$.
Computing the contribution of the remainder of the exceptional
locus to $P_{geom}$, as is done in Section \ref{toric}, allows us
to derive the contribution of this point, hence the motivic
integral $D^{(e_i)}$. However, we prefer a direct computation, in
order to obtain a shorter and more elementary method to compute
the geometric Poincar\'e series of a toric surface singularity.

 First, we take care of arcs $\psi$ for which
$\tilde{\psi}\notin \tilde{X}$. Their contribution equals
$\LL^{-2}(\LL-1)$ times
\begin{eqnarray*}\frac{(\LL-1)\LL^{-(3e_i+2)s-3\nu_2-1}}{(1-\LL^{-(e_i+1)s-\nu_2})(1-\LL^{-(2e_i+1)s-2\nu_2-1})}
+\LL^{-e_is-\nu_2}\frac{1+\LL^{-(e_i+1)s-\nu_2}}{1-\LL^{-(2e_i+1)s-2\nu_2-1}}\,.
\end{eqnarray*}

 The integral over $(E_i'\cap E''_i)\setminus
(E_i\cup H''_i)$ is equal to $\LL^{-3}(\LL-1)^{2}$ times
$$(\LL-2+\frac{(\LL-1)\LL^{-s-1}}{1-\LL^{-s-1}})\frac{\LL^{-(2e_i+1)s-2\nu_2}}{1-\LL^{-(2e_i+1)s-2\nu_2}}
\frac{1+\LL^{-(e_i+1)s-\nu_2}}{1-\LL^{-(e_i+1)s-\nu_2}}\,.$$

 We blow up
$E'_i\cap E''_i$ to compute the contribution of $H''_i\cap
E'_i\cap E''_i$. Let $G$ be the exceptional divisor. We denote the
strict transforms of $\bar{E}_i$, $E'_i$ and $E''_i$ again by the
same symbols. The question is how
\begin{eqnarray*}
\lefteqn{\max\,\{\min\{\lfloor
c(\bar{\psi},E'_i)/2+c(\bar{\psi},E''_i)/2\rfloor,
c(\bar{\psi},E'_i), c(\bar{\psi},E''_i)+c(\bar{\psi},H''_i)\},}
\\&&\qquad\qquad\qquad\qquad\qquad\qquad\qquad\qquad c(\bar{\psi},E'_i),\,c(\bar{\psi},E'_i)-c(\bar{\psi},H''_i)\,\}
\end{eqnarray*} behaves on $G$.

Straightforward computation shows that the term associated to
$H''_i\cap E'_i\cap E''_i$ is equal to
$$\LL^{-3}(\LL-1)^{2}\frac{\LL^{-(2e_i+1)s-2\nu_2}}{1-\LL^{-(2e_i+1)s-2\nu_2}}$$
times
\begin{eqnarray*}
\lefteqn{1+\frac{\LL^{-(e_i+1)s-\nu_2}}{1-\LL^{-(e_i+1)s-\nu_2}}
+\frac{(\LL-1)\LL^{-(3e_i+2)s-3\nu_2-1}}{(1-\LL^{-(e_i+1)s-\nu_2})(1-\LL^{-(2e_i+1)s-2\nu_2-1})}}
\\&&\qquad\qquad\qquad\qquad\qquad\qquad\qquad+\LL^{-e_is-\nu_2}\frac{1+\LL^{-(e_i+1)s-\nu_2}}{1-\LL^{-(2e_i+1)s-2\nu_2-1}}\,.
\end{eqnarray*}

 Finally, we compute the contribution of $E_i\cap E'_i\cap
E''_i$. Our expression for $d_{X}(\psi)$ reduces to
$$e_ic(\bar{\psi},E_i)+e_ic(\bar{\psi},E'_i)+e_ic(\bar{\psi},E''_i)+\max\{c(\bar{\psi},E'_i),c(\bar{\psi},E''_i)\}\,.$$
Hence, we get
$$\LL^{-3}(\LL-1)^{3}
\frac{\LL^{-(3e_i+1)s-2\nu_2-\nu_1}}{(1-\LL^{-(2e_i+1)s-2\nu_2})(1-\LL^{-e_is-\nu_1})}
\frac{1+\LL^{-(e_i+1)s-\nu_2}}{1-\LL^{-(e_i+1)s-\nu_2}}\,.$$

 Bringing all these terms
together, we see that
$$D^{(e_i)}(X,s)=\LL^{-3}(\LL-1)\frac{1+\LL^{-(e_i+1)s-\nu_2}}{1-\LL^{-(2e_i+1)s-2\nu_2}}
\{\frac{(\LL-1)^{2}\LL^{-(2e_i+1)s-2\nu_2}}{(1-\LL^{-s-1})(1-\LL^{-e_is-\nu_1})}+\LL^{-e_is-\nu_1}\}\,.$$
The contribution of $x_i$ to $P_{geom}$, i.e.
$$-\frac{\LL^{3}D^{(e_i)}(\LL^{3}T)}{1-\LL^{3}T}\,,$$
 is equal to
$$-\frac{(\LL-1)(1+\LL^{e_i+1}T^{e_i+1})}{(1-\LL^{3}T)(1-\LL^{2e_i-1}T^{2e_i+1})}\{\frac{(\LL-1)^{2}\LL^{2e_i-1}
T^{2e_i+1}}{(1-\LL^{2}T)(1-\LL^{e_i-1}T^{e_i})}+\LL^{e_i-1}T^{e_i}\}\,.$$

Now suppose $d>3$.
Consider the singular surface $U$ in $\A^{3}_k$, defined by
$yz-x^{2e_i+1}=0$. We consider $\A^{d}_k$ as the product of
$\A^{3}_k$ and $\A^{d-3}_k$, and we identify $\A^{3}_k$ with
$\A^{3}_k\times 0$. The minimal toric resolution of $U$ is induced
by a succession $g':A'\rightarrow \A^{3}_{k}$ of blow-ups of
points in $\A^{3}_{k}$, see Section \ref{toric}. Blowing up the
same points in $\A^{d}_k$ yields a morphism $g:A\rightarrow
\A^{d}_k$. The geometric Poincar\'e series $\tilde{P}_{U,geom}(T)$
of $U$ at the origin $O$ is defined intrinsically, hence does not
depend on the embedding in ambient space. It is computed in
Section \ref{toric}, only using our formula for $d=3$.

Applying the change of variables formula to the proper birational
morphism $g$, we write $D_{O}(U,s)$ as $D^{(e_i)}(U,s)+R(s)$,
where $D^{(e_i)}(U,s)$ is, of course, the integral over the arcs
in the resolution space with origin at the unique point where the
strict transform of $U$ and the exceptional locus intersect
non-transversally. The term $R$ is the motivic integral over the
remainder of the fiber of $g$ over the origin.

 It follows immediately from our formula that
$D^{(e_i)}(U,s)=D^{(e_i)}(X,s)$. The formula in Lemma 1 yields
$$P_{U,geom}(T)
=\frac{1-\LL^{d}(D^{(e_i)}(\LL^{d}T)+R(\LL^{d}T))}{1-\LL^{d}T}\,.$$
The term $\LL^{d}R(\LL^{d}T)$ can be easily computed, as is done
for $d=3$ in Section \ref{toric}. Writing $e_i$ as $2m_i+1$, and
$d$ as $t+2$, we see that
$$\LL^{d}D^{(e_i)}(\LL^{d}T)$$ is equal to
\begin{eqnarray*}
\lefteqn{-(1-\LL^{t+2}T)\frac{(\LL-1)\sum_{j=1}^{m_i-1}\LL^{j}T^{j}}{(1-\LL^{2}T)(1-\LL^{m_j-1}T^{m_j})}+(1-\LL^{t+2}T)
\frac{(\LL-1)\LL T}{(1-T)(1-\LL^{2}T)}}
\\&&+\frac{\LL^{2}(\LL^{t}-1)(\LL-1)T^{2}}{(1-T)(1-\LL^{2}T)}
-\frac{(\LL-1)^{2}([\Pro^{t}]-2)
\sum_{i=1}^{m_j-1}\LL^{i}T^{i+1}}{(1-T)(1-\LL^{m_j-1}T^{m_j})}
\\&&-2\frac{(\LL-1)^{2}(\LL^{t}-1)\LL^{2}T
\sum_{i=1}^{m_j-1}\LL^{i}T^{i+1}}{(1-T)(1-\LL^{m_j-1}T^{m_j})(1-\LL^{2}T)}\,.
\end{eqnarray*}
\section{Quasirational singularities}
The theorems in the previous section allow us to give a partial
answer to the question raised by Lejeune-Jalabert and
Reguera-Lopez in \cite{LejReg}.
\begin{lemma}
Let $C$ be an irreducible curve over $k$. If the class $[C]$ of
$C$ in the Grothendieck ring belongs to $\Z[\LL]$, then $C$ is
rational.
\end{lemma}
\begin{proof}
 It follows from Riemann-Roch that a smooth projective
curve of genus 0 is rational \cite{Hart}. The Poincar\'e
polynomial $P(u)$ is an additive invariant, i.e. it is
well-defined on $K_{0}(Var_k)$, since it is obtained from the
Hodge polynomial by identifying the two variables. For a complete
smooth curve $X$, the polynomial $P[X]$ is equal to
$u^{2}-2g(X)u+1$, where $g(X)$ is the genus of $X$.

A curve $C$ is rational if and only if its projective smooth
birational model $\bar{C}$ is. Furthermore, $[C]\in \Z[\LL]$
implies $[\bar{C}]\in \Z[\LL]$. It
 follows from the identities $P(\LL)=u^{2}$ and $P(1)=1$ that the linear term of
 $P[\bar{C}]$ is zero, hence $\bar{C}$ is rational, and so is $C$.
\end{proof}
\begin{theorem}\label{rat}
Let $x$ be an isolated singularity of a surface $X\subset Y$, with
$Y$ smooth, and assume that there exists an embedded resolution of
the germ of $X$ at $x$, satisfying the conditions of Theorem
\ref{formula} or \ref{form2}. The local Poincar\'e series of $X$
at $x$ can be written as a rational function with numerator and
denominator in $\Z[\LL][T]$, if and only if $x$ is a quasirational
 singularity.
\end{theorem}
Let us recall that a surface singularity $(X,x)$ is quasirational
if "only rational curves can come out of the singularity, no
matter how we blow it up birationally" \cite{Ab}. Rational surface
singularities are quasirational.
\begin{proof}
One implication is straightforward: suppose that $x$ is a
quasirational singularity. It suffices to write $P_{geom}$ in
terms of the motivic integral $D(s)$, and to observe that the
quasirationality of $x$ implies that all Grothendieck brackets of
the strata of the exceptional locus, emerging in the expression
for $P_{geom}$, belong to $\Z[\LL]$.

 Let $e_i$ be the minimal depth
of a global exceptional component $E_i$ in the resolution of
$(X,x)$, at which a non-rational exceptional component appears. We
will prove that the coefficient $A_i$ of $T^{e_i}$ in $P_{geom}$
is not contained in $\Z[\LL]$.

It is clear from our formulae that $A_i$ is equal to the sum of a
term in $\Z[\LL]$ with $(\LL-1)\sum \LL^{\mu_i} [C_i]$, where the
$\mu_i$ are positive integers, and we take the sum over all
non-rational exceptional components $C_i$ on $X'$ which are
contained in an exceptional divisor $E_i$ of depth $e_i$. Since
the coefficient of the linear term of $ P[C_i]$ is strictly
greater than zero, this sum can not be contained in $\Z[\LL]$, for
its Poincar\'e polynomial will contain a term of odd degree.
\end{proof}

\section{Toric surfaces}\label{toric}
Theorem \ref{formula} provides an elementary method to compute the
geometric Poincar\'e series of a toric surface singularity, which
is substantially shorter than the techniques used in
\cite{LejReg}.

 Throughout this section, $X$ is a singular
affine toric surface, defined by a cone $\sigma$ generated by
$(1,0)$ and $(p,q)$, where $0<p<q$ and $p,q$ are relatively prime.
Let $(b_1,\ldots,b_s)$ be the entries occurring in the
Hirzebruch-Jung continued fraction associated to $q/(q-p)$, and
$(c_1,\ldots,c_t)$ the components of the continued fraction of
$q/p$ \cite{F}\cite{Oda}. The relation between the $b_i$ and the
$c_j$ is explained in \cite{Nic1}. Let furthermore $\Theta$ be the
union of compact faces of the convex hull of $\sigma\cap
N\setminus 0$, and $\check{\Theta}$ be the union of compact faces
of the convex hull of $\check{\sigma}\cap M\setminus 0$.

The minimal resolution of $X$ is a toric modification induced by a
subdivision of $\sigma$ into simple cones. The vectors occurring
in this subdivision can be listed as follows:
$$v_0=(1,0),\,v_1=(1,1),\ldots,\,v_{j+1}=b_jv_{j}-v_{j-1},\ldots,\,v_{s+1}=b_sv_{s}-v_{s-1}=(p,q)\,.$$
The exceptional divisors $D_j\cong \Pro^{1}$ of this resolution
 correspond to the newly introduced vectors $v_j$, $j=1,\ldots,s$,
and $D_j$ is known to have self-intersection number $-b_j$.

The $c_j$ have a geometric significance of their own: subdividing
$\check{\sigma}$ into simple cones, i.e. taking the minimal set of
generators for the semi-group $\check{\sigma}\cap M$, yields an
embedding of $X$ into affine $(t+2)$-space; the ideal of $X$ is
generated by $x_{i-1}x_{i+1}-x_i^{c_i}$, $i=1,\ldots,t$.

We will factor the canonical toric resolution into a sequence of
blow-ups of zero-dimensional orbits, which can be immediately
extended to an embedded resolution for $X$ using the embedding in
affine space mentioned above. Blowing up the unique
zero-dimensional orbit $O$ of $V$ corresponds, by \cite{Kempf}, to
the toric modification corresponding to the subdivision $\Sigma$
of $\sigma$ introducing all  primitive vectors normal to the edges
of $\check{\Theta}$. This comes down to inserting $v_1$,
$v_{s-1}$, and all $v_i$ determining vertices of $\Theta$, i.e.
the $v_i$ for which $b_i\neq 2$. Let $a$ be the number of vectors
introduced in $\Sigma$, i.e. the number of elements in
$\{b_2,\ldots,b_{s-1}\}$ differing from 2 augmented by two, and
let $b=a-r-1$ be the number of pairs of adjacent vectors in
$\Sigma$, that is, pairs of vectors in $\Sigma$ with multiplicity
$1$. The number $b$ is equal to the number of $c_j$ equal to 3,
while $r$ is equal to the cardinality of $\{c_j>3\}$.

The singularities left after blowing up $O$ are all rational
singularities of type $A_{c}$. In fact, they are recovered from
the $b_i$ by omitting $b_1$ and $b_s$, and isolating all sequences
of 2's in the remaining $b_i$. Let $c$ be the number of 2's in
such a sequence. This number $c$ can be recovered from the $c_j$:
it is equal to $c_j-3$, with $j$ chosen such that the vertex of
$\check{\Theta}$ corresponding to $x_j$ lies on the two edges
whose normal directions determine the cone in our fan $\Sigma$
corresponding to this sequence of 2's. Moreover, each of the $c_j$
which is bigger than 3 will induce a singularity in this way. The
singularity will be resolved after blowing up the zero-dimensional
orbit corresponding to the associated singular cone (thus
inserting 2 vectors, or 1 if $c=1$) and repeating this procedure
$\lfloor c/2 \rfloor$ times. If $c$ is even, we get a chain of
exceptional divisors intersecting $\tilde{X}$ transversally; if
$d_k$ is odd, we get an intersection point of multiplicity 2 in
the last stage of the resolution process.

This factorization allows us to embed our resolution in ambient
affine space, simply by blowing up the corresponding points in
this space. Let $h:\tilde{Y}\rightarrow Y=\A_{k}^{t+2}$ be the
proper birational morphism obtained in this way, and let
$\tilde{X}$ be the strict transform of $X$ (thus $\tilde{X}$ is
the canonical resolution surface). The points of $\tilde{X}$ where
there's no transversal intersection with the exceptional locus of
$h$ correspond to adjacent vectors in the simple subdivision of
$\sigma$ which are introduced in one and the same blow-up.

In \cite{Nic1}, we used this embedded resolution to compute the
motivic Igusa Poincar\'e series of $X$. Our computation of
$P_{geom}$ will be very similar.

 Define $E_{-1}$ to be the strict transform
of $X$ under $h$. Let $E_0$ be the strict transform of the
exceptional divisor that is created in the first blow-up, and let
$E_{i,j}$ be the strict transform of the exceptional divisor
induced by the $j$-th blow-up of the singularity corresponding to
the $i$-th sequence of 2's in $b_2,\ldots,b_{s-1}$.

 Let $c'_i$ be the $i$-th component of $(c_1,\ldots,c_t)$ which is strictly larger than 3, and put $d_i$ equal to $c'_i-3$. We let $I$ denote the index set
$$\{-1,0\}
\cup \{(i,j)\,|\,i\in \{1,\ldots,r\},\,j\in\{1,\ldots,\lceil
d_i/2\rceil \}\}\,.$$  We stratify $\tilde{X}$ in the usual way:
for each subset $J$ of $I$, we define $E_{J}$ to be
$\cap_{\alpha\in J}E_{\alpha}$, while $E_{J}^{o}$ denotes
$E_{J}\setminus \cup_{\alpha\notin J}E_{\alpha}$.

We attach to each $E_{\alpha}$ a pair of numerical data
$(N_{\alpha},\nu_{\alpha})$ as follows: $$(N_{-1},\nu_{-1})=(1,t),
\,(N_{0},\nu_{0})=(1,t+2), 
(N_{(i,j)},\nu_{(i,j)})=(j+1,(j+1)(t+1)+1)\,.$$
 Then
\begin{eqnarray*}
D(s)&=&\LL^{-(t+2)}\sum_{J\subset
I,J\nsubseteq\{-1\}}[E_{J}^{o}]\prod_{\alpha\in J}
\frac{(\LL^{codim\,E_{\alpha}}-1)\LL^{-{N_{\alpha}
s-\nu_{\alpha}}}}{1 -\LL^{-{N_{\alpha} s-\nu_{\alpha}}}}
\\&& \qquad
+b\,\LL^{-(t+2)}\{D^{(1)}(s)-\frac{(\LL^{t}-1)(\LL-1)\LL^{-3s-2t-2}}{(1-\LL^{-s-t})(1-\LL^{-2s-t-2})}\}
\\&&+\sum_{d_i\
even}\{D^{(d_i/2+1)}(s)-\frac{(\LL^{t}-1)(\LL-1)\LL^{-(d_i/2+2)s-(d_i/2+2)(t+1)}}{(1-\LL^{-s-t})(1-\LL^{-(d_i/2+1)s-(d_i/2+1)(t+1)-1})}\}
\end{eqnarray*}
The last terms in the expression for $D(s)$ correct for
non-transversal intersection. We refer to \cite{Nic1} for a more
detailed description of the terms.

We immediately recover a result from \cite{LejReg}, stating that
$P_{geom}$ is trivial when all $c_j$ are equal to 2, i.e. that
$P_{geom}$ equals the local geometric series of a smooth point.
More generally, we can state the following corollary of Theorem
\ref{formula}:
\begin{corollary}
Let $X\subset Y$ be varieties over $k$, where $Y$ is smooth, and
$X$ has dimension $m$. Let $x$ be an isolated singularity of $X$,
and let $h:Y'\rightarrow Y$ be a composition of blow-ups of
points, satisfying the conditions of Theorem \ref{formula}, with
exceptional divisor $E$, which is a linear chain of components
$E_i$. Let $X'$ be the strict transform of $X$. Assume that, for
each $i$, the class of $E_i\cap X'$ in $K_{0}(Var_k)$ is equal to
$[\Pro^{m-1}_{k}]$, and $[E_i\cap E_{i+1}\cap X]=[\Pro^{m-2}_k]$.
 Then $P_{geom}$ is trivial.
\end{corollary}
To put it intuitively: if the embedded resolution looks like $X$
was smooth at $x$ all along, $P_{geom}$ cannot distinguish $x$
from a smooth point. The general idea is that similar embedded
resolution graphs yield similar geometrical Poincar\'e series.

 Let us reduce the formula for
 $P_{geom}(T)$
to the expression given in \cite{LejReg}. First, we treat the case
$t=1$, that is, we compute the geometric Poincar\'e series
associated to the toric surface singularity defined by
$xz=y^{c_j}$. We already know what happens for $c_j=2$, so we may
suppose $c_j>2$. Let $D'(s)$ be the sum of terms of $D(s)$
corresponding to index sets not containing any of the couples
$(i,j)$. The contribution
$$\frac{1-\LL^{d}D'(\LL^{d}T)}{1-\LL^{d}T}$$
is equal to
\begin{eqnarray*}
\frac{1}{1-\LL^{2}T}+\frac{(\LL-1)\LL
T}{(1-T)(1-\LL^{2}T)}+\frac{\LL^{2}(\LL-1)^{2}T^{2}}{(1-T)(1-\LL^{2}T)(1-\LL^{3}T)}\,.
\end{eqnarray*}
If $c_j=3$, we have to add the contribution
$$-\frac{(\LL-1)(1+\LL^{2}T^{2})}{(1-\LL^{3}T)(1-\LL T^{3})}\{\frac{(\LL-1)^{2}\LL
T^{3}}{(1-\LL^{2}T)(1-T)}+T\}\,,$$ and we obtain
$$P_{geom}(T)=\frac{1}{1-\LL^{2}T}+(\LL-1)\frac{(\LL-1)T+\LL T^{2}}{(1-\LL^{2}T)(1-\LL
T^{3})}\,.$$

If $c_j>3$, write $c_j$ as $2m_j+n_j$, where $m_j$ is a positive
integer, and $n_j$ is either zero or one.

First, we consider the strata $[E_{(i,j)}^{o}]=(\LL-1)^{2}$ and
$[E^{o}_{\{0,(1,j)\}}]=[E^{o}_{\{(i,j),(i+1,j)\}}]=\LL-1$. The
corresponding terms in the expression for $P_{geom}$ amount to
$$-\frac{(\LL-1)^{3}}{1-\LL^{3}T}\{\frac{\LL T^{3}}{(1-T)(1-\LL T^{2})}
+\sum_{i=2}^{m_j-1}\frac{\LL^{i-1}T^{i}}{(1-\LL^{i-1}T^{i})(1-\LL^{i}T^{i+1})}\}-\frac{(*)^{(j)}}{1-\LL^{3}T}$$
where the term $(*)_j$ depends on the value of $n_j$. The part
between braces is easily seen to be equal to
$$\frac{\LL^{m_j-1}T^{m_j+1}+\sum_{i=1}^{m_j-2}\LL^{i}T^{i+1}}{(1-T)(1-\LL^{m_j-1}T^{m_j})}\,.$$

If $n_j=0$, the term $(*)^{(j)}$ equals
$$(\LL-1)(\LL^{2}-\LL+1)\frac{\LL^{m_j-1}T^{m_j}}{1-\LL^{m_j-1}T^{m_j}}\,,$$
while, if $n_j=1$, it is equal to
$$(\LL-1)^{3}\frac{\LL^{m_j-1}T^{m_j}}{1-\LL^{m_j-1}T^{m_j}}\,.$$
Hence, we get
$$-\frac{(\LL-1)^{3}\sum_{i=1}^{m_j-1}\LL^{i}T^{i+1}}{(1-\LL^{3}T)(1-T)(1-\LL^{m_j-1}T^{m_j})}+(n_j-1)\frac{\LL(\LL-1)\LL^{m_j-1}T^{m_j}}{(1-\LL^{3}T)(1-\LL^{m_j-1}T^{m_j})}\,.$$

Now, we look at the remaining strata, except for the points where
we get non-transversal intersection. This yields, in exactly the
same way,
\begin{eqnarray*}\lefteqn{-2\frac{(\LL-1)^{3}\LL^{2}T}{(1-\LL^{2}T)(1-\LL^{3}T)}
\frac{\sum_{i=1}^{m_j-1}\LL^{i}T^{i+1}}{(1-T)(1-\LL^{m_j-1}T^{m_j})}}
\\&&\quad\quad\qquad -(n_j-1)\frac{(\LL-1)^{3}\LL^{m_j+1}T^{m_j+1}}{(1-\LL^{m_j-1}T^{m_j})(1-\LL^{2}T)(1-\LL^{3}T)}\,.
\end{eqnarray*}
If $c_j$ is even, we see that $P_{geom}(T)$ equals
$$\frac{1}{1-\LL^{2}T}+(\LL-1)\frac{\sum_{i=1}^{m_j-1}\LL^{i}T^{i}}{(1-\LL^{2}T)(1-\LL^{m_j-1}T^{m_j})}\,.$$
If $c_j$ is odd, we have to include the contribution of the point
where there's no transversal intersection with the exceptional
locus. This yields
$$P_{geom}(T)=\frac{1}{1-\LL^{2}T}+(\LL-1)\frac{(1+\LL^{m_j-1}T^{m_j})\sum_{i=1}^{m_j}\LL^{i}T^{i}-\LL^{m_j-1}T^{m_j}}{(1-\LL^{2}T)(1-\LL^{2m_j-1}T^{2m_j+1})}\,.$$

To conclude, let us consider the case $d>3$. It follows from
Theorems \ref{formula} and \ref{form2}, our computations for
$d=3$, and the fact that the geometric Poincar\'e series is
defined intrinsically, that the contribution to $P_{geom}$ of all
strata, except for $E_{\{0\}}^{o}$ and $E_{\{-1,0\}}^{o}$, is
equal to
\begin{eqnarray*}
\lefteqn{-(a-1)\frac{(\LL-1)\LL
T}{(1-T)(1-\LL^{2}T)}-(a-1)\frac{\LL^{2}(\LL^{t}-1)(\LL-1)T^{2}}{(1-T)(1-\LL^{2}T)(1-\LL^{t+2}T)}}
\\&&+\sum_{c_j>2 \
even}(\LL-1)\frac{\sum_{i=1}^{m_j-1}\LL^{i}T^{i}}{(1-\LL^{2}T)(1-\LL^{m_j-1}T^{m_j})}
\\&&+\sum_{c_j>2\ odd}(\LL-1)\frac{(1+\LL^{m_j-1}T^{m_j})\sum_{i=1}^{m_j}\LL^{i}T^{i}-\LL^{m_j-1}T^{m_j}}{(1-\LL^{2}T)(1-\LL^{2m_j-1}T^{2m_j+1})}\,.
\end{eqnarray*}
This observation allows us to conclude that
\begin{eqnarray*}
P_{geom}(T)&=&\frac{1}{1-\LL^{2}T} +\sum_{c_j>2 \
even}(\LL-1)\frac{\sum_{i=1}^{m_j-1}\LL^{i}T^{i}}{(1-\LL^{2}T)(1-\LL^{m_j-1}T^{m_j})}
\\&&+\sum_{c_j>2\ odd}(\LL-1)\frac{(1+\LL^{m_j-1}T^{m_j})\sum_{i=1}^{m_j}\LL^{i}T^{i}-\LL^{m_j-1}T^{m_j}}{(1-\LL^{2}T)(1-\LL^{2m_j-1}T^{2m_j+1})}\,.
\end{eqnarray*}

 It follows from
results in \cite{DLinvent} that this formula holds not only over
$\hat{\mathcal{M}}_k$, but already over $\mathcal{M}_{k}$.
 Of course, this
formula also holds for surface singularities with "the same
embedded resolution graph" as a toric surface singularity.

\section{The arithmetic series}
In this section, $k$ is an arbitrary field of characteristic zero,
not necessarily algebraically closed, and we denote by $k^{alg}$
an algebraic closure. Let $X$ be a variety over $k$, and let $x$
be a point of $X(k)$. As we've seen, the local geometric
Poincar\'e series counts the jets in $\mathcal{L}_{n}(X)_{x}$,
which can be lifted to an arc in $\mathcal{L}(X)_{x}$. However,
working scheme-theoretically, we allow extensions of our base
field $k$ in this lifting process (which is necessary to ensure
that $[\pi_{n}\mathcal{L}(X)_{x}]$ is well-defined). Hence,
$P_{geom}(T)$ is insensitive to issues of rationality. These are
taken into account by the arithmetic series $P_{arith}(T)$.

Bittner \cite{Bitt} gave a short proof of the existence of a ring
morphism $\chi_{mot}$ from
 the Grothendieck ring $K_{0}(Var_k)$ of varieties over $k$ to the Grothendieck
 ring $K_{0}(CH_k)$ of Chow motives over $k$, sending the class of a smooth
 projective variety to the class of its associated Chow motive,
 and sending $\LL$ to the class of the Tate motive $\LL_{mot}$.
 The existence of this map was originally proven in \cite{GuiNav}.
In \cite{DL1}, Denef and Loeser constructed a morphism
$$\chi_{c}:K_{0}(PFF_k)\rightarrow K_{0}^{mot}(Var_k)\otimes \Q\
,$$ where $K_{0}^{mot}(Var_k)$ is the image of $K_{0}(Var_k)$
under the morphism $\chi_{mot}$, and $K_{0}(PFF_k)$ is the
Grothendieck group of pseudo-finite fields containing $k$.
Elements of $K_{0}(PFF_k)$ are equivalence classes of ring
formulas over $k$.

For this construction, it is important to understand the structure
of $K_{0}(PFF_k)$. The theory of quantifier elimination for
pseudo-finite fields \cite{FriJar}\cite{FriSac}, states that
quantifiers can be eliminated if one adds some relations to the
language, which have a geometric interpretation in terms of Galois
covers. This interpretation yields a construction for $\chi_{c}$.
It is important for our purposes that, if our original ring
formula $\varphi$ did not contain any quantifiers in the first
place, $\chi_{c}$ maps $[\varphi]$ to the class of the
constructible set defined by $\varphi$ in $K^{mot}_{0}(Var_k)$. We
refer to \cite{Nic1} for a short introduction to the arithmetic
series, and to \cite{DL}\cite{DL1} for the original work on
arithmetic integration.

Let us define the local arithmetic Poincar\'e series $P_{arith}$
of $X$ at $x$. Since we're working locally, we may assume that $X$
is a subvariety of some affine space $Y=\A^{d}_k$. It follows from
Greenberg's theorem that we can find, for each positive integer
$n$, a ring formula $\varphi_{n}$ over $k$, such that, for all
fields $K$ containing $k$, the $K$-rational points of
$\mathcal{L}_{n}(X)_{x}$ that can be lifted to a $K$-rational
point of $\mathcal{L}(X)_{x}$, correspond to the tuples satisfying
the interpretation of $\varphi_{n}$ in $K$. We define the
arithmetic Poincar\'e series to be
$$P_{arith}(T)=\sum_{n\geq 0}\chi_{c}([\varphi_{n}])T^{n}\ .$$ As was proven in \cite{DL}, it is
rational over $K_{0}^{mot}(Var_k)[\LL^{-1}]\otimes\Q$.

Suppose that there exists an embedded resolution
$h:\tilde{Y}\rightarrow Y$ for $X$, defined over $k$, which
satisfies, after a base change to $k^{alg}$, the conditions of
Theorem \ref{formula} or \ref{form2}. We demand that each
component of the exceptional locus contains a $k$-rational point
of the strict transform $\tilde{X}$ that does not lie on any other
exceptional component, and that the divisors $F_1$ and $F_2$ in
Theorem \ref{form2} are defined over $k$. The main result of this
section is the following theorem.
\begin{theorem}\label{arith}
If these conditions are satisfied, the motivic series $P_{arith}$
and $P_{geom}$ are equal in $K_{0}^{mot}(Var_k)\otimes \Q[[T]]$.
\end{theorem}
\begin{proof}
We will prove that, for each field $K$ containing $k$, and each
positive integer $n$, a $K$-rational point $j_n$ of
$\mathcal{L}_n(X)_{x}$ lifts to an arc on $X$ if and only if it
lifts to a $K$-rational arc $\psi$ on $X$. Hence, $\varphi_n$ and
the set of quantifier-free equalities and inequalities describing
the constructible set $\pi_n\mathcal{L}(X)_x$ define the same
class in $K_0(PFF_k)$, which proves Theorem \ref{arith}.

So suppose $j_n$ lifts to an arc on $X$, and let $\eta$ be a
$K$-rational arc on $Y$, lifting $j_n$. By definition, $d_X(\eta)>
n$. Theorems \ref{formula} and \ref{form2} not only give a formula
for $d_X(\eta)$, but also give you an optimal approximation for
$\eta$ in $\mathcal{L}(X)$. Now it suffices to observe that this
approximation can be chosen to be $K$-rational.
\end{proof}

We recover a particular case of a theorem in \cite{Nic1}\,:
\begin{corollary}
If $(X,x)$ is the germ of a toric surface singularity, $P_{arith}$
and $P_{geom}$ are equal in $K_{0}^{mot}(Var_k)\otimes \Q[[T]]$.
\end{corollary}
\bibliographystyle{hplain}
\bibliography{wanbib,wanbib2}
\end{document}